\documentclass[12pt]{article}%
\usepackage{graphicx}
\usepackage[intlimits]{amsmath}
\usepackage{latexsym}
\usepackage{amsfonts}
\usepackage{amssymb}%
\makeatletter
\newfont{\footsc}{cmcsc10 at 8truept}
\newfont{\footbf}{cmbx10 at 8truept}
\newfont{\footrm}{cmr10 at 10truept}
\newtheorem{theorem}{Theorem}

\newtheorem{claim}[theorem]{Claim}

\newtheorem{conjecture}[theorem]{Conjecture}

\newenvironment{proof}[1][Proof]{\noindent{\textbf {#1}  }}  {\hfill$\Box$\bigskip}

\def\blfootnote{\xdef\@thefnmark{}\@footnotetext}
\begin{document}

\title{A new type of Ramsey-Tur\'{a}n problems}
\author{H. Li\thanks{{\small LRI, Universit\'{e}.de Paris-sud and CNRS, Orsay
F-91405}} ,\ V. Nikiforov\thanks{{\small Department of Mathematical Sciences,
University of Memphis, Memphis, TN 38152}} \thanks{Research has been supported
by NSF Grant \# DMS-0906634.} \ and R.H. Schelp\footnotemark[2]\\{\small e-mail: li@lri.fr; vnikifrv@memphis.edu; rschelp@memphis.edu}}
\maketitle

\begin{abstract}
We introduce and study a new type of Ramsey-Tur\'{a}n problems, a typical
example of which is the following one:

Let $\varepsilon>0$ and $G$ be a graph of sufficiently large order $n$ with
minimum degree $\delta\left(  G\right)  >3n/4.$ If the edges of $G$ are
colored in blue or red, then for all $k\in\left[  4,\left\lfloor \left(
1/8-\varepsilon\right)  n\right\rfloor \right]  ,$ there exists a
monochromatic cycle of length $k.$\medskip


\end{abstract}

\subsection*{Introduction and main results}

In notation we follow \cite{Bol98}. Given a graph $G$ with edge set $E\left(
G\right)  ,$ a\emph{ }$2$\emph{-coloring} of $G$ is a partition $E\left(
G\right)  =E\left(  R\right)  \cup E\left(  B\right)  ,$ where $R$ and $B$ are
subgraphs of $G$ with $V\left(  R\right)  =V\left(  B\right)  =V\left(
G\right)  .$

It is well-known (see \cite{FaSc74},\cite{Ros73}) that if $n>3$ and the edges
of the complete graph $K_{2n-1}$ are colored in two colors, there is a
monochromatic cycle of length $k$ for every $k\in\left[  3,n\right]  .$ This
is best possible since the complete bipartite graph $K_{n-1,n-1}$ and its
complement contain no $n$-cycle when $n$ is odd; yet the statement can be
improved further: as proved in \cite{NiSc08}, essentially the same conclusion
follows from weaker premises:\medskip

\emph{Let }$G$\emph{ be a graph of sufficiently large order }$2n-1,$
\emph{with }$\delta\left(  G\right)  \geq\left(  2-10^{-6}\right)  n.$\emph{
If }$E\left(  G\right)  =E\left(  R\right)  \cup E\left(  B\right)  $\emph{ is
a }$2$\emph{-coloring,} \emph{then }$C_{k}\subset R$\emph{ for all }%
$k\in\left[  3,n\right]  $\emph{ or }$C_{k}\subset B$\emph{ for all }%
$k\in\left[  3,n\right]  .$\medskip

This assertion suggests a new class of Ramsey-Tur\'{a}n extremal problems,
which somewhat do not fit the traditional framework of this area as given,
e.g., in \cite{SiSo01}. In this note, we will be mainly concerned with the
following conjecture.

\begin{conjecture}
\label{con1}Let $n\geq4$ and let $G$ be a graph of order $n$ with
$\delta\left(  G\right)  >3n/4.$ If $E\left(  G\right)  =E\left(  R\right)
\cup E\left(  B\right)  $ is a $2$-coloring, then $C_{k}\subset R$ or
$C_{k}\subset B$ for all $k\in\left[  4,\left\lceil n/2\right\rceil \right]
.$
\end{conjecture}

If true, this conjecture is tight as shown by the following example: let
$n=4p,$ color the edges of the complete bipartite graph $K_{2p,2p}$ in blue,
and insert a red $K_{p,p}$ in each of its vertex classes. The union of the
blue and red edges gives a graph $G$ with $\delta\left(  G\right)  =3n/4,$ but
clearly the coloring produces no monochromatic odd cycle. In fact, in the
concluding remarks we describe a more sophisticated coloring scheme that gives
$2^{n^{2}/8-O\left(  n\log n\right)  }$ non-isomorphic colorings of $G$ with
no monochromatic odd cycle.

Note that under the premises of Conjecture \ref{con1} we don't necessarily get
a monochromatic triangle: indeed, it is long known (see \cite{Sos69}) that
splitting $K_{5}$ into a blue and a red $5$-cycle and then blowing-up this
coloring, we obtain a graph $G$ with $\delta\left(  G\right)  =4n/5$ with no
monochromatic triangle.

The abundance of nonisomorphic extremal colorings suggests that the complete
proof of Conjecture \ref{con1} might be difficult. In this note we give a
partial solution of the conjecture, stated in the theorem below.

\begin{theorem}
\label{th1}Let $\varepsilon>0,$ let $G$ be a graph of sufficient large order
$n,$ with $\delta\left(  G\right)  >3n/4.$ If $E\left(  G\right)  =E\left(
R\right)  \cup E\left(  B\right)  $ is a $2$-coloring, then $C_{k}\subset R$
or $C_{k}\subset B$ for all $k\in\left[  4,\left\lfloor \left(
1/8-\varepsilon\right)  n\right\rfloor \right]  .$
\end{theorem}

We have to admit that even the proof of this weaker statement needs quite a
bit of work. Before proving the theorem we introduce some notation and a
number of supporting results.

\subsubsection*{Some graph notation}

Given a graph $G,$ we write:

- $G\left[  U\right]  $ for the graph induced by a set $U\subset V\left(
G\right)  ;$

- $\Gamma_{G}\left(  u\right)  $ for the set of neighbors of a vertex $u\in
V\left(  G\right)  $ and $d_{G}\left(  u\right)  $ for $\left\vert \Gamma
_{G}\left(  u\right)  \right\vert ;$

- $\Delta\left(  G\right)  $ for the maximum degree of $G;$

- $ec\left(  G\right)  $ and $oc\left(  G\right)  $ for the lengths of the
longest even and odd cycles in $G.$

\subsubsection*{\textbf{Preliminary results}}

Here we state several known results that will be needed throughout our
proof.\medskip

\textbf{Theorem A (Bondy \cite{Bon71}) }\emph{If }$G$\emph{ is a graph of
order }$n$\emph{, with }$\delta\left(  G\right)  >n/2,$\emph{ then }$G$\emph{
contains }$C_{k}$\emph{ for all }$k\in\left[  3,n\right]  $\emph{.}\medskip

In \cite{Bol78}, p. 150, Bollob\'{a}s gave the following size version of
Theorem A.

\textbf{Theorem B (Bollob\'{a}s \cite{Bol78}) }\emph{If }$G$\emph{ is a graph
of order }$n$\emph{,} \emph{with }$e\left(  G\right)  >n^{2}/4,$\emph{ then
}$G$\emph{ contains }$C_{k}$\emph{ for all }$k\in\left[  3,\left\lceil
n/2\right\rceil \right]  $\emph{.}\medskip

\textbf{Theorem C (H\"{a}ggkvist \cite{Hag82}) }\emph{Let }$k\geq2$\emph{ be a
n integer and let }$G$\emph{ be a nonbipartite, }$2$\emph{-connected graph of
sufficiently large order }$n.$\emph{ If }%
\[
\delta\left(  G\right)  >\frac{2}{2k+1}n,
\]
\emph{then }$G$\emph{ contains a }$C_{2k-1}$\emph{.}\medskip

\textbf{Theorem D (Gould, Haxell, Scott \cite{GHS02}) }For all $c>0,$ there
exists $K=K\left(  c\right)  $ such that if $G$ is a graph of order $n>K$,
with $\delta\left(  G\right)  \geq cn,$ then $C_{k}\subset G$ for all even
$k\in\left[  4,ec\left(  G\right)  -K\right]  $ and all odd $k\in\left[
K,oc\left(  G\right)  -K\right]  .$\medskip

\textbf{Theorem E (Erd\H{o}s, Gallai \cite{ErGa59}) }Let $k\geq1.$ If
$e\left(  G\right)  >k\left\vert G\right\vert /2$, then $G$ contains a path
and a cycle of length at least $k+1.$\medskip

\subsubsection*{Proof of Theorem \ref{th1}}

Let $G$ be a graph of sufficiently large order $n$ with minimum degree
$\delta>3n/4,$ and let $V=V\left(  G\right)  $ be the vertex set of $G.$
Suppose $E\left(  G\right)  =E\left(  R\right)  \cup E\left(  B\right)  $ is a
$2$-coloring, and let $\varepsilon>0.$ Clearly it is enough to prove the
theorem for $\varepsilon$ sufficiently small, so we shall assume that
$\varepsilon$ is as small as needed.

Assume that $n$ is sufficiently large and that there is no monochromatic
$k$-cycle for some $k\in\left[  4,\left\lfloor \left(  1/8-\varepsilon\right)
n\right\rfloor \right]  $. For convenience we shall refer to this assumption
as \emph{the main assumption.} In the following claims, we establish a number
of properties of $R$ and $B$ that follow from the main assumption.

\begin{claim}
\label{cl1}Both graphs $R$ and $B$ are nonbipartite.
\end{claim}

\begin{proof}
Indeed, suppose that, say, $B$ is bipartite and let $U$ and $W$ be its vertex
classes. By symmetry, assume that $\left\vert U\right\vert \geq\left\vert
W\right\vert $ and let $R_{1}=R\left[  U\right]  .$ Then for every $u\in U,$%
\[
d_{R_{1}}\left(  u\right)  \geq\delta\left(  G\right)  -\left\vert
W\right\vert >\frac{3}{4}n-n+\left\vert U\right\vert \geq\frac{\left\vert
U\right\vert }{2}.
\]
Now, Theorem A implies that $R_{1}$ contains $C_{k}$ for all $k\in\left[
3,\left\vert U\right\vert \right]  .$ Since $\left\vert U\right\vert \geq
n/2,$ this inequality contradicts the main assumption, completing the proof of
Claim \ref{cl1}.
\end{proof}

\begin{claim}
\label{cl2}$e\left(  R\right)  >n^{2}/8$ and $e\left(  B\right)  >n^{2}/8$.
\end{claim}

\begin{proof}
Assume, by symmetry, that $e\left(  R\right)  \leq n^{2}/8.$ Then
\[
e\left(  B\right)  \geq\frac{1}{2}\delta n-e\left(  R\right)  >\frac{3}%
{8}n^{2}-\frac{1}{8}n^{2}=\frac{1}{4}n^{2}.
\]
Theorem B now implies that $B$ contains a $k$-cycle for all $k\in\left[
3,\left\lceil n/2\right\rceil \right]  ,$ contradicting the main assumption
and proving the claim.
\end{proof}

\begin{claim}
\label{cl3} $\delta\left(  R\right)  >2\varepsilon n$ and $\delta\left(
B\right)  >2\varepsilon n.$
\end{claim}

\begin{proof}
Assume, by symmetry, that there is a vertex $u$ such that $d_{R}\left(
u\right)  \leq2\varepsilon n.$ Then
\begin{equation}
d_{B}\left(  u\right)  =\delta\left(  G\right)  -d_{R}\left(  u\right)
>\left(  \frac{3}{4}-2\varepsilon\right)  n. \label{i1}%
\end{equation}
Let $U=\Gamma_{B}\left(  u\right)  $. Clearly $\Gamma_{G}\left(  v\right)
\cap U=\Gamma_{G}\left(  v\right)  \backslash\left(  V\backslash U\right)  ,$
and so%
\[
d_{R\left[  U\right]  }\left(  v\right)  +d_{B\left[  U\right]  }\left(
v\right)  =\left\vert \Gamma_{G}\left(  v\right)  \cap U\right\vert
\geq\left\vert \Gamma_{G}\left(  v\right)  \right\vert -\left\vert V\backslash
U\right\vert >\frac{3}{4}n-n+\left\vert U\right\vert =\left\vert U\right\vert
-\frac{1}{4}n.
\]
Therefore,
\[
e\left(  R\left[  U\right]  \right)  +e\left(  B\left[  U\right]  \right)
>\frac{1}{2}\left(  \left\vert U\right\vert -\frac{1}{4}n\right)  \left\vert
U\right\vert .
\]

If $e\left(  R\left[  U\right]  \right)  >\left\vert U\right\vert ^{2}/4,$
then by Theorem B, $R$ contains $C_{k}$ for all $k\in\left[  3,\left\lceil
\left\vert U\right\vert /2\right\rceil \right]  ,$ and for sufficiently small
$\varepsilon$, this contradicts the main assumption. Thus $e\left(  R\left[
U\right]  \right)  \leq\left\vert U\right\vert ^{2}/4.$ On the other hand, if
$B\left[  U\right]  $ contains a path $P$ of order $\left\lfloor \left(
1/8-\varepsilon\right)  n\right\rfloor ,$ we again reach a contradiction,
since $P$, together with the vertex $u,$ gives a $k$-cycle for all
$k\in\left[  3,\left\lfloor \left(  1/8-\varepsilon\right)  n\right\rfloor
\right]  .$ Therefore, Theorem E implies that $e\left(  R\left[  U\right]
\right)  \leq\left(  1/8-\varepsilon\right)  n\left\vert U\right\vert /2.$ In
summary,%
\[
\frac{1}{2}\left(  \left\vert U\right\vert -\frac{1}{4}n\right)  \left\vert
U\right\vert <e\left(  R\left[  U\right]  \right)  +e\left(  B\left[
U\right]  \right)  \leq\frac{1}{4}\left\vert U\right\vert ^{2}+\frac{1}%
{2}\left(  \frac{1}{8}-\varepsilon\right)  n\left\vert U\right\vert ,
\]
implying that
\[
\left\vert U\right\vert -\frac{1}{4}n<\frac{1}{2}\left\vert U\right\vert
+\left(  \frac{1}{8}-\varepsilon\right)  n
\]
and so,
\[
\left\vert U\right\vert <\left(  \frac{3}{4}-2\varepsilon\right)  n,
\]
a contradiction with (\ref{i1}), completing the proof of the claim.
\end{proof}

\begin{claim}
\label{cl4} At least one of the graphs $R$ or $B$ is $2$-connected.
\end{claim}

\begin{proof}
Assume for a contradiction that both $R$ and $B$ are at most $1$-connected.
That is to say, we can remove a vertex $v$ such that $R-v$ is disconnected and
a vertex $u$ such that $B-u$ is disconnected. Note that the components of
$R-v$ and $B-u$ are at least of size $\varepsilon n$, as follows from the
previous claim. Letting $S=\left\{  u,v\right\}  $ and $W=V\backslash S,$ we
see that we can remove a set $S$ of at most two vertices so that the graphs
$R^{\prime}=R\left[  W\right]  $ and $B^{\prime}=B\left[  W\right]  $ are
disconnected. Hence $R^{\prime}=R_{1}\cup R_{2}$ and $B^{\prime}=B_{1}\cup
B_{2},$ where $R_{1}$ and $R_{2}$ are vertex disjoint graphs, and so are
$B_{1}$ and $B_{2}.$ Set
\[
V_{1}=R_{1}\cap B_{1},\text{ }V_{2}=R_{1}\cap B_{2},\text{ }V_{3}=R_{2}\cap
B_{1},\text{ }V_{4}=R_{2}\cap B_{2}.
\]
Clearly, $W=\cup_{i=1}^{4}V_{i}$ is a partition of $W.$ Note also that there
are no cross edges in $G$ between $V_{1}$ and $V_{4},$ and also between
$V_{2}$ and $V_{3}.$ We thus have
\[
d_{G\left[  W\right]  }\left(  w\right)  \leq\left\{
\begin{array}
[c]{cc}%
n-\left\vert S\right\vert -1-\left\vert V_{4}\right\vert  & \text{if }w\in
V_{1}\\
n-\left\vert S\right\vert -1-\left\vert V_{1}\right\vert  & \text{if }w\in
V_{4}\\
n-\left\vert S\right\vert -1-\left\vert V_{2}\right\vert  & \text{if }w\in
V_{3}\\
n-\left\vert S\right\vert -1-\left\vert V_{3}\right\vert  & \text{if }w\in
V_{2}%
\end{array}
\right.
\]
This immediately implies that
\[
3n-3\left\vert S\right\vert -4\geq4\delta\left(  G\left[  W\right]  \right)
\geq4\left(  \delta\left(  G\right)  -\left\vert S\right\vert \right)
>3n-4\left\vert S\right\vert ,
\]
a contradiction, completing the proof of the claim.
\end{proof}

\subsubsection*{Existence of monochromatic short odd cycles}

Below we shall prove that, under the premises of the theorem, there is a
monochromatic $\left(  2k+1\right)  $-cycle for every $k\geq2,$ provided
$n>n_{0}\left(  k\right)  $. Unfortunately this cannot give the full proof
since $n_{0}\left(  k\right)  $ grows faster than linear in $k$. In the claim
below and the subsequent short argument, we dispose of all odd cycles that are
longer than $5.$

\begin{claim}
\label{cl5}If $\Delta\left(  B\right)  \geq n/2+4k,$ then $B$ contains a
$C_{2k+1}.$
\end{claim}

\begin{proof}
Suppose that $\Delta\left(  B\right)  \geq n/2+4k,$ choose a vertex $u$ with
$d_{B}\left(  u\right)  =\Delta\left(  B\right)  ,$ and write $U$ for the set
of its neighbors. If $B\left[  U\right]  $ contains a path of order $2k,$ this
path, together with $u,$ gives a cycle $C_{2k+1},$ so the claim is proved in
this case. On the other hand, $e\left(  R\left[  U\right]  \right)
\leq\left\vert U\right\vert ^{2}/4$, for otherwise Theorem B implies that
$R\left[  U\right]  $ has cycles of all lengths from $3$ to $\left\vert
U\right\vert /2,$ and since $\left\vert U\right\vert /2\geq n/4+2k>n/8,$ this
contradicts the main assumption. Hence,
\[
e\left(  R\left[  U\right]  \right)  +e\left(  B\left[  U\right]  \right)
\leq\frac{1}{4}\left\vert U\right\vert ^{2}+k\left\vert U\right\vert ,
\]
On the other hand, as in the proof of Claim \ref{cl3}, we see that
\[
e\left(  R\left[  U\right]  \right)  +e\left(  B\left[  U\right]  \right)
\geq\frac{1}{2}\left(  \delta\left(  G\right)  -\left(  n-\left\vert
U\right\vert \right)  \right)  \left\vert U\right\vert >\frac{1}{2}\left(
\left\vert U\right\vert -\frac{1}{4}n\right)  \left\vert U\right\vert ,
\]
and so
\[
\frac{1}{2}\left(  \left\vert U\right\vert -\frac{1}{4}n\right)  \left\vert
U\right\vert <\frac{1}{4}\left\vert U\right\vert ^{2}+k\left\vert U\right\vert
,
\]
implying in turn $\left\vert U\right\vert <n/2+4k,$ contrary to our selection
of $U.$ This completes the proof.
\end{proof}

In view of Claim \ref{cl5}, to prove the existence of monochromatic $\left(
2k+1\right)  $-cycle we can assume that $\Delta\left(  B\right)  <n/2+4k,$ and
consequently, $\delta\left(  R\right)  >n/4-4k.$ Note that $n/4-4k>2n/\left(
2k+1\right)  $ for every $k\geq4$ and $n$ sufficiently large. Then, since $R$
is $2$-connected, Theorem C implies that $R$ contains $C_{2k-1}$ for every
$k\geq4$ provided $n$ is sufficiently large.

It remains to prove that there is a monochromatic $C_{5};$ this is the bulk of
our effort. We shall assume that $\delta\left(  R\right)  >n/4-16,$ and for a
contradiction let us assume that there is no monochromatic $C_{5}.$ The proof
is split in three cases.\bigskip

\emph{(i) }$R$ \emph{contains no triangle}\bigskip\textbf{ }

Theorem C implies that there is a $7$-cycle $C$ in $R,$ and since
$C_{3}\nsubseteq R$ and $C_{5}\nsubseteq R$, the cycle $C$ must be induced.
Let $v_{1},\ldots,v_{7}$ be its vertices listed consecutively along the cycle,
and note that $v_{1},v_{3},v_{5},v_{7},v_{2},v_{4},v_{6}$ is also a $7$-cycle
in $B;$ we shall refer to it by $C^{\prime}.$ Note that no vertex of\ $V$ can
be joined in $R$ to three vertices of $C,$ for otherwise we have either
$C_{3}\subset R$ or $C_{5}\subset R$. Hence%
\[
\sum_{i=1}^{7}d_{R}\left(  v_{i}\right)  \leq2n.
\]
That is to say,
\[
\sum_{i=1}^{7}d_{B}\left(  v_{i}\right)  >7\cdot\frac{3n}{4}-\sum_{i=1}%
^{7}d_{R}\left(  v_{i}\right)  \geq\left(  7\cdot\frac{3}{4}-2\right)
n\geq\frac{13n}{4}.
\]
Therefore, there is a vertex in $V$, joined to $4$ vertices of $C^{\prime}$ in
$B$. An easy check shows that $C_{5}\subset B$, a contradiction. Thus, $R$
must contain a triangle.\bigskip

\emph{(ii) }$R$ \emph{contains triangles, but no two triangles share an
edge}\bigskip

Let $v_{1},v_{2},v_{3}$ be the vertices of a triangle, and $U_{1},U_{2},U_{3}$
be their neighborhoods in $R.$ Clearly our premise implies that
\[
\left\vert U_{1}\cap U_{2}\right\vert =\left\vert U_{2}\cap U_{3}\right\vert
=\left\vert U_{1}\cap U_{3}\right\vert =1.
\]
Also letting
\[
U_{1}^{\prime}=U_{1}\backslash\left\{  u_{2},u_{3}\right\}  ,\text{ }%
U_{2}^{\prime}=U_{2}\backslash\left\{  u_{1},u_{3}\right\}  ,\text{ }%
U_{3}^{\prime}=U_{3}\backslash\left\{  u_{1},u_{2}\right\}  ,
\]
we see that there are no red cross edges between $U_{i}^{\prime}$ and
$U_{j}^{\prime},$ $\left(  1\leq i<j\leq3\right)  ,$ for otherwise
$C_{5}\subset R.$

Since $C_{5}\nsubseteq R,$ none of the graphs $R\left[  U_{i}^{\prime}\right]
$ contains a path of length $3;$ thus, Theorem E implies that $e\left(
R\left[  U_{i}^{\prime}\right]  \right)  \leq\left\vert U_{i}^{\prime
}\right\vert .$ That is to say, for every $i=1,2,3,$ we can choose a vertex
$v_{i}\in U_{i}^{\prime}$ such that
\[
\left\vert \Gamma_{R}\left(  v_{i}\right)  \backslash U_{i}^{\prime
}\right\vert \geq d_{R}\left(  v_{i}\right)  -\delta\left(  R\left[
U_{i}^{\prime}\right]  \right)  >d_{R}\left(  v_{i}\right)  -2>n/4-18.
\]
Now, for every $i=1,2,3$, set $W_{i}=\Gamma_{R}\left(  v_{i}\right)
\backslash U_{i}^{\prime},$ and note that $W_{i}\cap U_{j}^{\prime
}=\varnothing$ and $W_{i}\cap W_{j}=\varnothing$ for all $1\leq i<j\leq3.$
Indeed, $W_{i}\cap U_{j}=\varnothing$ as there are no red cross-edges between
$U_{i}^{\prime}$ and $U_{j}^{\prime}.$ Likewise if, say $w\in W_{1}\cap
W_{2},$ then $w,v_{1},u_{1},u_{2},v_{2}$ is a $5$-cycle since $w\neq
u_{1},u_{2}.$

Now we see that the $6$ sets $U_{i}^{\prime},W_{i},$ $\left(  i=1,2,3\right)
$ are pairwise disjoint and so,
\[
n\geq\sum\left\vert U_{i}^{\prime}\right\vert +\left\vert W_{i}\right\vert
\geq6\left(  n/4-18\right)  ,
\]
which is a contradiction for $n$ large enough, completing the proof in this
case. It remains to consider the last possibility:\bigskip

\emph{(iii) }$R$ \emph{contains two triangles sharing an edge}\bigskip

Let $v_{1},v_{2},v_{4}$ and $v_{2},v_{3},v_{4}$ be the vertices of these
triangles, and let and $U_{1},U_{2},U_{3}$ be the neighborhoods of
$v_{1},v_{2},v_{3}$ in $R.$ Clearly since $C_{5}\nsubseteq R$, then
\[
\left\vert U_{1}\cap U_{2}\right\vert \leq2,\text{ }\left\vert U_{2}\cap
U_{3}\right\vert \leq2,\text{ }\left\vert U_{1}\cap U_{3}\right\vert =2.
\]
Also letting
\[
U_{1}^{\prime}=U_{1}\backslash\left\{  u_{2},u_{4},u_{3}\right\}  ,\text{
}U_{2}^{\prime}=U_{2}\backslash\left\{  u_{1},u_{3},u_{4}\right\}  ,\text{
}U_{3}^{\prime}=U_{3}\backslash\left\{  u_{1},u_{2},u_{4}\right\}  ,
\]
we see that there are no cross edges in $R$ between $U_{i}^{\prime}$ and
$U_{j}^{\prime},$ $\left(  1\leq i<j\leq3\right)  ,$ for otherwise
$C_{5}\subset R.$

Again, since $C_{5}\nsubseteq R,$ none of the graphs $R\left[  U_{i}^{\prime
}\right]  $ contains a path of length $3;$ thus, Theorem E implies that
$e\left(  R\left[  U_{i}^{\prime}\right]  \right)  \leq\left\vert
U_{i}^{\prime}\right\vert .$ That is to say, for every $i=1,2,3,$ we can
choose a vertex $v_{i}\in U_{i}^{\prime}$ such that
\[
\left\vert \Gamma_{R}\left(  v_{i}\right)  \backslash U_{i}^{\prime
}\right\vert \geq d_{R}\left(  v_{i}\right)  -\delta\left(  R\left[
U_{i}^{\prime}\right]  \right)  \geq d_{R}\left(  v_{i}\right)  -2>n/4-18.
\]
Now, for every $i=1,2,3$, set $W_{i}=\Gamma_{R}\left(  v_{i}\right)
\backslash U_{i}^{\prime},$ and note that $W_{i}\cap U_{j}^{\prime
}=\varnothing$ and $W_{i}\cap W_{j}=\varnothing$ for all $1\leq i<j\leq3.$
Indeed, $W_{i}\cap U_{j}=\varnothing,$ as there are no red cross-edges between
$U_{i}^{\prime}$ and $U_{j}^{\prime}.$ Likewise if, say $w\in W_{1}\cap
W_{2},$ then $w,v_{1},u_{1},u_{2},v_{2}$ is a $5$-cycle since $w\neq
u_{1},u_{2}.$

Now we see that the $6$ sets $U_{i}^{\prime},W_{i},$ $\left(  i=1,2,3\right)
$ are pairwise disjoint and so
\[
n\geq\sum\left\vert U_{i}^{\prime}\right\vert +\left\vert W_{i}\right\vert
\geq3\left(  n/4-19\right)  +3\left(  n/4-18\right)  ,
\]
which is a contradiction for $n$ large enough, completing the proof of the
existence of monochromatic short odd cycles.

\subsubsection*{Existence of all cycles}

At this point we know that $R$ is a $2$-connected nonbipartite graph, with
$e\left(  R\right)  >n^{2}/8$ and $\delta\left(  R\right)  \geq2\varepsilon
n.$ Let $K=K\left(  2\varepsilon\right)  ,$ where $K\left(  \cdot\right)  $ is
the function from Theorem D and set $L=2\left\lceil 1/2\varepsilon\right\rceil
+K$. Based on the above properties of $R$ we shall prove that $C_{k}\subset R$
for all even $k\in\left[  4,n/8-L\right]  $ and all odd $k\in\left[
K,n/8-L\right]  .$ Our main tool will be Theorem D, but to apply it we have to
prove the following claim.

\begin{claim}
\label{cl6}$ec\left(  R\right)  >n/8-2k$ and $oc\left(  R\right)  >n/8-2k.$
\end{claim}

\begin{proof}
In view of
\[
\delta\left(  R\right)  >2\varepsilon n=\frac{2}{2/\left(  2\varepsilon
\right)  }n>\frac{2}{2\left(  \left\lceil 1/2\varepsilon\right\rceil \right)
+1}n,
\]
Theorem C implies that $R$ contains a $\left(  2r-1\right)  $-cycle
$C^{\prime}$, where $r=\left\lceil 1/2\varepsilon\right\rceil ;$ let say
$v_{1},\ldots,v_{2r-1}$ be the vertices of $C^{\prime}$ listed consecutively
along the cycle. Write $R^{\prime}$ for the graph obtained from $R$ by
omitting the set $\left\{  v_{1},\ldots,v_{2r-1}\right\}  .$ We have
\[
e\left(  R^{\prime}\right)  >e\left(  R\right)  -\left(  2r-1\right)  \left(
n-1\right)  >\frac{n^{2}}{8}-\left(  2r-1\right)  \left(  n-1\right)
\geq\left(  \frac{n}{8}-2r\right)  \left(  n-2r+1\right)  ,
\]
the last inequality holding for $n$ sufficiently large. Now, Theorem E implies
that $R^{\prime}$ contains an $l$-cycle $C^{\prime\prime}$ for some
$l>n/4-4r;$ let say $u_{1},\ldots,u_{l}$ be its vertices listed consecutively
along the cycle. Since $R$ is $2$-connected, there exist two vertex disjoint
paths $P^{\prime}$ and $P^{\prime\prime}$ joining two vertices $v_{i},v_{j}$
of $C^{\prime}$ to two vertices $u_{s},u_{t}$ of $C^{\prime\prime}$ and having
no other vertices in common with either $C^{\prime}$ or $C^{\prime\prime}.$
The vertices $u_{s},u_{t}$ split $C^{\prime\prime}$ into two paths with common
ends $u_{s}$ and $u_{t}$. Taking the longer of these paths, together with
$P^{\prime}$ and $P^{\prime\prime},$ we obtain a path $P$ joining $v_{i}%
,v_{j}$ and having no other vertices in common with $C^{\prime}$. Note that
the length of $P$ is at least $\left\lceil l/2\right\rceil +2>n/8-2r.$ Finally
note that $v_{i}$ and $v_{j}$ split $C^{\prime}$ into two paths joining
$v_{i}$ and $v_{j}$ - one of even length and one of odd length; these paths,
together with $P,$ give one odd and one even cycle, each longer than $n/8-2r.$
\end{proof}

Applying Theorem D with $c=2\varepsilon,$ and letting $K=K\left(
2\varepsilon\right)  ,$ we see that $C_{k}\subset R$ for all even $k\in\left[
4,ec\left(  G\right)  -K\right]  $ and all odd $k\in\left[  K,oc\left(
G\right)  -K\right]  .$ In view of Claim \ref{cl6}, we see that $C_{k}\subset
R$ for all even $k\in\left[  4,n/8-L\right]  $ and all odd $k\in\left[
K,n/8-L\right]  ,$ as stated.

Now, to complete the proof of Theorem \ref{th1}, recall that for all $k\geq3$
and $n$ sufficiently large, either $C_{2k-1}\subset R$ or $C_{2k-1}\subset B$.
Therefore, for all $k\in\left[  4,n/8-L\right]  ,$ either $C_{k}\subset R$ or
$C_{k}\subset B,$ a contradiction with the main assumption, completing the
proof of the theorem.$\hfill\square$

\subsection*{Concluding remarks}

Let $U_{1},U_{2},U_{3},U_{4}$ be sets of size $p.$ Let $G$ be the graph with
vertex set $\cup_{i=1}^{4}U_{i},$ and let the edges of $G$ be all $U_{i}%
-U_{j}$ edges for $i\neq j.$ Setting $n=4p,$ we shall show that there exist
$2^{n^{2}/8-O\left(  n\log n\right)  }$ non-isomorphic edge colorings of $G$
which do not produce any monochromatic cycle. Indeed, color all $U_{1}-U_{2}$
and $U_{3}-U_{4}$ edges in blue; color all $U_{1}-U_{3}$ and $U_{2}-U_{4}$
edges in red; color all $U_{1}-U_{4}$ and $U_{2}-U_{3}$ edges arbitrarily in
blue or red. It is easy to check that both the red and blue graphs are
bipartite, and so there is no monochromatic odd cycle. Since the $U_{1}-U_{4}$
and $U_{2}-U_{3}$ edges can be colored in $2^{2p^{2}}$ ways, we obtain at
least $2^{n^{2}/8}/n!=2^{n^{2}/8-O\left(  n\log n\right)  }$ different
colorings of $G$ all of which avoid a monochromatic odd cycle.

We conclude with another conjecture, which seems a bit easier than Conjecture
\ref{con1}.

\begin{conjecture}
Let $0<c<1$ and $G$ be a graph of sufficiently large order $n.$ If
$\delta\left(  G\right)  >cn,$ then for every $2$-coloring of $E\left(
G\right)  $, there is a monochromatic $C_{k}$ for some $k\geq cn$.
\end{conjecture}

\end{document}